# Relationship formula between nonlinear polynomial equations and the corresponding Jacobian matrix

W. Chen

**Present mail address** (as a JSPS Postdoctoral Research Fellow): Apt.4, West 1st floor, Himawari-so, 316-2, Wakasato-kitaichi, Nagano-city, Nagano-ken, 380-0926, JAPAN

Permanent affiliation and mail address: Dr. Wen CHEN, P. O. Box 2-19-201, Jiangshu University of Science & Technology, Zhenjiang City, Jiangsu Province 212013, P. R. China

Present e-mail: chenw@homer.shinshu-u.ac.jp

Permanent email: chenwwhy@hotmail.com

**Abstract:** This paper provides a general proof of a relationship theorem between nonlinear analogue polynomial equations and the corresponding Jacobian matrix, presented recently by the present author. This theorem is also verified generally effective for all nonlinear polynomial algebraic system of equations. As two particular applications of this theorem, we gave a Newton formula without requiring the evaluation of nonlinear function vector as well as a simple formula to estimate the relative error of the approximate Jacobian matrix. Finally, some possible applications of this theorem in nonlinear system analysis are discussed.

**1. Introduction**:

Recently, the present author [1] proved a relationship theorem as stated below

**Theorem 1.** If $N^{(m)}(U)$ and $J^{(m)}(U)$ are defined as nonlinear numerical analogue of the m order nonlinear differential operator and its corresponding Jacobian matrix, respectively, then $N^{(m)}(U) = \frac{1}{m} J^{(m)}(U)U$ is always satisfied irrespective if which numerical technique is employed to discretize.

Some significant applications of the above theorem were given in [1]. However, that paper failed to provide a general proof of this theorem including fractional order nonlinear problems. The objective of this paper is to verify the effectiveness of this theorem for general nonlinear polynomial function vector. The theorem was also used to derive a Newton formula without requiring the evaluation of nonlinear function vector as well as a simple formula to estimate the relative error of the approximate Jacobian matrix.

## 2. General proof

The method of weighted residuals (MWR) is recognized the origin of almost all popular numerical techniques [2, 3]. Consider the differential equations of the form

$$\psi\{u\} - f = 0, \text{ in } \Omega \tag{1}$$

with the following boundary conditions

$$u = \bar{u}, \quad \text{on } \Gamma_u \tag{2a}$$

$$q = \bar{q}, \quad \text{on } \Gamma_q \tag{2b}$$

where n is the outward normal to the boundary, $\Gamma = \Gamma_u + \Gamma_q$, the upper bars indicate known values on boundary, and $q = \partial u / \partial n$. More complex boundary conditions can be easily included but they will not be considered here for the sake of brevity. In the MWR, the desired function u in the differential governing equation is first approximated by a set of linearly independent basis functions $\phi_k(x)$, such that

$$u = \hat{u} = \sum_{k=1}^{n} c_k \phi_k, \tag{3}$$

where $c_k$'s are the unknown parameters. In the Galerkin and FE methods, the basis functions are usually chosen so as to satisfy certain given conditions such as the boundary conditions and the continuity. In addition, these basis functions should be complete.

Substituting Eq. (3) into Eq. (1) produces an error, which is called the residual, namely,

$$\psi\{\hat{u}\} - f = R \neq 0. \tag{4}$$

This error or residual R is forced to be zero in the average sense by setting weighted integral of the residuals equal to zero, namely,

$$\int_\Omega [\psi\{\hat{u}\} - f] W_j d\Omega = \int_\Omega R W_j d\Omega = 0, \quad j=1, 2, ..., N, \tag{5}$$

where $W_j$'s are weighting functions. The use of different weighting functions and basis functions give rise to different numerical techniques such that the Galerkin, Least square, finite element, boundary element, moments, spectral methods, finite difference and collocations methods.

This paper places its emphasis on the nonlinear computations. Let us consider (1+s)-order nonlinear operator of general form:

$$p(u)[r(u)]^s + L(u) = f, \tag{6}$$

where p(u), r(u) and L(u) are linear differential operators, f denotes the constant, and s any real number. The scheme of weighted residuals approximations Eq. (6) is given by

$$\int_\Omega \left[ p(\hat{u})[r(\hat{u})]^s + L(\hat{u}) - f \right] W_i d\Omega = \int_{\Gamma_2} (\hat{q} - \bar{q}) W_i d\Gamma - \int_{\Gamma_1} (\hat{u} - \bar{u}) \frac{\partial W_i}{\partial n} d\Gamma, \quad j=1,2,...,N. \tag{7}$$

Substituting equation (3) into the above equation (7) yields

$$\psi(c) = \int_\Omega \left[ \left( \sum_{k=1}^n c_k p(\phi_k) \right) \left( \sum_{k=1}^n c_k r(\phi_k) \right)^s + \sum_{k=1}^n c_k L(\phi_k) - f \right] W_i d\Omega - $$
$$\int_{\Gamma_2} \left( \sum_{k=1}^n c_k \frac{\partial \phi_k}{\partial u} - \bar{q} \right) W_i d\Gamma + \int_{\Gamma_1} \left( \sum_{k=1}^n c_k \phi_k - \bar{u} \right) \frac{\partial W_i}{\partial n} d\Gamma = 0 \quad j=1, 2, ..., N. \tag{8}$$

The above equation (8) can be restated as

$$\psi(c) = Dc + N^{(1+s)}(c) + b = 0, \tag{9}$$

where c is a vector comprised of the unknown $c_k$, D is the linear coefficient matrix, and $N^{(1+s)}(c)$ means the (1+s)-order nonlinear vector term. For example, it is the quadratic nonlinear function vector when s=1 and 3/2 order nonlinear one when s=1/2. b represents the contact vector. The Jacobian matrix of nonlinear vector $N^{(1+s)}(c)$ is given by

$$J_{ij}^{(1+s)} = \frac{\partial N^{(1+s)}(c)}{\partial c_j} = \int_\Omega \left[ p(\phi_j)\left(\sum_{k=1}^n c_k r(\phi_k)\right)^s + sr(\phi_j)\left(\sum_{k=1}^n c_k r(\phi_k)\right)^{s-1}\left(\sum_{k=1}^n c_k p(\phi_k)\right) \right] W_i d\Omega \quad i,j=1,\ldots,N. \quad (10)$$

By using equations (8), (9) and (10), we can easily establish

$$N^{(1+s)}(c) = \frac{1}{1+s} J^{(1+s)}(c)c. \quad (11)$$

The above equation (11) is generally effective for numerical analogue of both integer and fractional order nonlinear differential or integral operators.

As mentioned earlier, nearly all popular numerical techniques such as the finite element, boundary element, finite difference, Galerkin, spectral, least square, moments and collocation methods and their variants can be derived from the weighted residual methodology. The only difference among a variety of numerical methods lies in the use of different weighting and basis functions. From the foregoing deduction, it is noted that equation (11) can be obtained no matter what weighting and basis functions we use in the method of weighted residuals. Therefore, it is straightforward that the formula (11) is generally effective for all numerical methods which can be derived from the weighted residual method.

In what follows, we will provide another straightforward proof of theorem 1 for nonlinear algebraic equations with explicit expression. At first, the concept of the Hadamard product and power are defined as in [4, 5].

**Definition 2.1.** Let matrices $A=[a_{ij}]$ and $B=[b_{ij}] \in C^{N \times M}$, the Hadamard product of matrices is defined as $A \circ B = [a_{ij} b_{ij}] \in C^{N \times M}$, where $C^{N \times M}$ denotes the set of $N \times M$ real matrices.

**Definition 2.2.** If matrix $A=[a_{ij}] \in C^{N \times M}$, then $A^{\circ q}=[a_{ij}^q] \in C^{N \times M}$ is defined as the Hadamard power of matrix A, where q is a real number. Especially, if $a_{ij} \neq 0$, $A^{\circ(-1)}=[1/a_{ij}] \in C^{N \times M}$ is defined as the Hadamard inverse of matrix A. $A^{\circ 0}=11$ is defined as the Hadamard unit matrix in which all elements are equal to unity.

**Theorem 2.1**: letting A, B and $C \in C^{N \times M}$, then

1. $A \circ B = B \circ A$            (12a)
2. $k(A \circ B)=(kA) \circ B$, where k is a scalar.            (12b)
3. $(A+B) \circ C = A \circ C + B \circ C$            (12c)

Consider the nonlinear $p(u)[r(u)]^s$ in equation (6) again, the corresponding numerical analogue by using a point-wise approximation technique can be expressed

$$p(u)[r(u)]^s = \{p(u)\}_i \circ \{r(u)\}_i^s = (A_p U) \circ (A_r U)^{\circ s} = N^{(1+s)}(U), \quad i=1,2,\ldots,N, \qquad (13)$$

where i indexes the number of discrete points; $A_p$ and $A_r$ represent the coefficient matrices of operators p(u) and r(u), respectively, dependent on specified numerical discretization scheme. The Hadamard product is exploited to express nonlinear discretization term in the above equation. It is noted that we use the desired function value vector U here instead of the unknown parameter vector C in equation (11). In fact, both are equivalent. The above explicit matrix formulation (13) is obtained in a straightforward and intuitive way. For more details see [6, 7]. The point-wise approximating numerical techniques include the finite difference, finite volume, collocation methods and their variants such as differential quadrature and pseudo-spectral methods. In addition, the numerical techniques based on radial basis functions can also express their analogue of nonlinear differential operators in the Hadamard product form. On the other hand, it is worth stressing that all explicit nonlinear polynomial equations which may not be originated from the numerical approximation can be expressed in the Hadamard product form.

The SJT product was introduced by the present author [6, 7] to efficiently compute analytical solution of the Jacobian derivative matrix.

**Definition 2.3.** If matrix $A=[a_{ij}]\in C^{N\times M}$, vector $U=\{u_j\}\in C^{N\times 1}$, then $A\Diamond U=[a_{ij}u_i]\in C^{N\times M}$ is defined as the postmultiplying SJT product of matrix A and vector U, where $\Diamond$ represents the SJT product. If M=1, $A\Diamond B=A\circ B$.

**Definition 2.4.** If matrix $A=[a_{ij}] \in C^{N \times M}$, vector $V=\{v_j\} \in C^{M \times 1}$, then $V^T \lozenge A = [a_{ij} v_j] \in C^{N \times M}$ is defined as the SJT premultiplying product of matrix A and vector V.

Considering the Hadamard nonlinear expression (13), we have

$$J^{(1+s)}(U) = \frac{\partial}{\partial U}\left\{(A_p U) \circ (A_r U)^{\circ s}\right\} = A_p \lozenge (A_r U)^{\circ s} + s A_r \lozenge (A_r U)^{\circ (s-1)} \lozenge (A_p U). \tag{14}$$

Formula (14) produces the analytical Jacobian matrix through simple algebraic computations. The SJT premultiplying product is related to the Jacobian matrix of the nonlinear formulations such as

$$q(U^k) = AU^{\circ k} = N^{(k)}(U), \tag{15}$$

where q is a linear operator, and m is any real number, and A is the numerical coefficient matrix of operator q( ). The corresponding Jacobian matrix is given by

$$J^{(k)}(U) = \frac{\partial}{\partial U}\{AU^k\} = \left(kU^{\circ(k-1)}\right)^T \lozenge A. \tag{16}$$

It is observed from the above formulas (15) and (16) that the Jacobian matrix of the nonlinear algebraic equations of the Hadamard product and power expression can be calculated by using the SJT product in the chain rules similar to those in differentiation of a scalar function. The computational effort of a SJT product is only $n^2$ scalar multiplication.

The finite difference method is often employed to calculate the approximate solution of the Jacobian matrix and also requires $O(n^2)$ scalar multiplications. In fact, the SJT product approach requires $n^2$ and $5n^2$ less multiplication operations than one and two order finite differences, respectively. Moreover, the SJT product produces the analytic solution of the Jacobian matrix. In contrast, the approximate Jacobian matrix yielded by the finite difference method affects the accuracy and convergence rate of the Newton-Raphson method, especially for highly nonlinear problems. The efficiency and accuracy of the SJT product approach were numerically demonstrated in [6-8].

We notice the following fact that the SJT product is closely related with the ordinary product of matrices, namely, if matrix $A=[a_{ij}] \in C^{N \times M}$, vector $U=\{u_i\} \in C^{N \times 1}$, then the postmultiplying SJT product of matrix A and vector U satisfies

$$A \Diamond U = \text{diag}\{u_1, u_2, \ldots, u_N\} A, \tag{17}$$

where matrix $\text{diag}\{u_1, u_2, \ldots, u_N\} \in C^{N \times N}$ is a diagonal matrix. Similarly, for the SJT premultiplying product, we have

$$V^T \Diamond A = A \, \text{diag}\{v_1, v_2, \ldots, v_M\}, \tag{18}$$

where vector $V=\{v_j\} \in C^{M \times 1}$.

By using equations (13), (14), (17) and (18), we can easily constitute

$$N^{(1+s)}(U) = \frac{1}{1+s} J^{(1+s)}(U) U. \tag{19}$$

Similarly, by using equations (15), (16), (17) and (18), we have

$$N^{(k)}(U) = \frac{1}{k} J^{(k)}(U) U. \tag{20}$$

The Hadamard expression as well as the previous MWR approximation of nonlinear algebraic equations encompasses all implicit and explicit polynomial nonlinear systems of equations. According to equations (11), (19) and (20), we can generalize theorem 1 to any nonlinear polynomial function vector as stated below:

**Theorem 2.2.** If $N^{(m)}(U)$ and $J^{(m)}(U)$ are respectively defined as the m-order nonlinear polynomial function vector and its corresponding Jacobian matrix, then $N^{(m)}(U) = \frac{1}{m} J^{(m)}(U) U$ is always satisfied.

The nonlinear polynomial function vector here indicate the systems in which all unknown variables are included only in polynomial functions. As an example of theorem 2.2, let us consider the typical quadratic nonlinear term

$$\varphi(X) = XAX \tag{21}$$

often encountered in control engineering, where X and A are rectangular matrices. Its Jacobian matrix is given by [5]

$$J = \frac{\partial \varphi}{\partial X} = I \otimes (AX) + (AX)^T \otimes I. \tag{22}$$

It is easily validated

$$\varphi(X) = \frac{1}{2} J \bar{X}, \tag{23}$$

where $\bar{X}$ is a vector by stacking the rows of matrix X. The above result is in agreement with theorem 2.2.

## 3. Newton iterative formula without evaluation of nonlinear function vector

The Newton method is of vital importance in nonlinear algebraic computations of various numerical algebraic analogoue equations. The time-consuming effort includes the repeated evaluation of the nonlinear function vector, Jacobian matrix and its inverse [9]. It is recognized that the repeated numerical integration of force vectors, in other words, nonlinear function vector, is one of factors which influence the efficiency of some popular numerical techniques [10]. Although the computational burden in this aspect does not take a crucial part of nonlinear solutions, the relative cost is absolutely nontrivial in solution of a large nonlinear system especially for the finite element scheme. By using theorem 2.2, the present author [11] derived a Newton iteration formula for nonlinear equation (9)

$$c^{k+1} = \frac{s}{1+s}c^k - \frac{1}{1+s}J^{-1}(c^k)(sDc^k + (1+s)b), \tag{24}$$

where u's with superscript k+1 and k mean respectively the iteration solutions at the (k+1)-th and k-th steps. It is noted that the iterative formula (24) does not require the evaluation of the nonlinear function vector yet is in fact equivalent to the standard Newton iterative formula. The simplicity and effectiveness of the formula (24) was verified in the solution of the static Navier-Stokes equations of the driven cavity problem. The present iterative formula may be especially useful to improve efficiency of various modified Newton methods, in which the repeated evaluation of the Jacobian and its inverse has been reduced or avoided. The cost in the repeated calculation of nonlinear function vector occupies a larger part of the total computational effort [9].

## 4. Error estimator of approximate Jacobian

The repeated calculation of the Jacobian matrix is required in general nonlinear computations. Function differences is one of conventional numerical methods often used in practice for this task. However, the method suffers from possible inaccuracy, particularly if the problem is highly nonlinear. Therefore, it is practically important to detect the accuracy of the approximate Jaocbian matrix by the finite difference method. By means of theorem 2.2, we here give a simple error estimator of the Jacobian matrix . Consider nonlinear equation (6), we have

$$err = \|\overline{\psi}(c) - \hat{J}(c)c\| / \|\overline{\psi}(c)\|, \tag{25}$$

where $\hat{J}(U)$ is the approximate Jacobian matrix of ψ(c). $\overline{\psi}(c)$ is defined

$$\overline{\psi}(c) = Dc + (1+s)N^{(1+s)}(c) + b, \tag{26}$$

which is different from ψ(c) of equation (9) in that the nonlinear vector $N^{(1+s)}$(c) is multiplied by the nonlinear order number 1+s. The above formula (25) provides a practically significant approach to examine the relative deviation error between the approximate and exact Jacobian matrices by vector norm.

## 4. Some remarks

The known von Karman equations of geometrically nonlinear bending, buckling and vibration of plate and shell are a mixed quadratic and cubic nonlinear differential system of equations. The resulting numerical discretization of von Karman equations by any numerical techniques can be written as

$$Lu + N^{(2)}(u) + N^{(3)}(u) = f, \tag{27}$$

where u is the desired displacement vector, L is the linear coefficient matrix, $N^{(2)}(u)$ and $N^{(3)}(u)$ mean respectively the quadratic and cubic nonlinear vector terms, and f represents the force vector. By using theorem 2.2, we have

$$Lu + \frac{1}{2}J^{(2)}(u)u + \frac{1}{3}J^{(3)}(u)u = f, \tag{28}$$

where $J^{(2)}(u)$ and $J^{(3)}(u)$ denote the Jacobian matrices of $N^{(2)}(u)$ and $N^{(3)}(u)$, respectively. Therefore, we can write equation (27) in linear like form

$$K(u)u = f, \tag{29}$$

where K(u) is definite physical stiffness matrix of system rather the common tangent stiffness matrix

$$\overline{K} = L + J^{(2)}(u) + J^{(3)}(u) \tag{30}$$

resulting from a Newton linearization procedure. The present K(u) can reflect better the original physical characteristics of nonlinear system.

The above analysis shows that theorem 2.2 can bring a direct analysis of nonlinear system without using linearization procedure such as the Newton method. In fact, all complex problems in physics and engineering such as shock, chaotic and soliton wave involve nonlinear system of governing equations. By applying theorem 2.2 to numerical discretization of these equations, we can get explicit system matrix. It is expected that special structure features of system matrix have close relations to the system definite behavior. We suggest an instance (sample) analysis approach in which the analysis of system matrix is done at certain time and space instance. The standard linear algebraic analysis such eigenvalues, condition number and circulant structures, etc. can be used to detect such system matrix. This may lead to a sensible algebraic understanding of mathematical physics features of the nonlinear systems.

Theorem 2.2 was applied in [1] to stability analysis of nonlinear initial value problems and direct solution of nonlinear algebraic equations by using the linear iterative methods such as the Jacobi, Gauss-Seidel and SOR methods without the linearization procedure. In particular, it should be pointed out that by applying theorem 2.2, we can get the explicit system matrix of nonlinear equations and thus can analyze stability behavior of nonlinear initial value problems by the existing effective approach for linear varying-coefficient problems [12].

Finally, it is worth stressing that a very large class of real-world nonlinear problems can be modeled or numerically discretized polynomial algebraic system of equations. Theorem 2.2 is in general applicable for all these problems. Therefore, this work is of practical significance in broad physical and engineering areas.